\documentclass[11pt,a4paper]{article}

\title{\bf The Weitzenb\"ock formula on the Wiener space and its application to
the asymptotic estimate of entropy}
\author{Dejun Luo\footnote{Email: luodj@amss.ac.cn}
\vspace{3mm}\\
{\footnotesize Key Lab of Random Complex Structures and Data
Science, Academy of Mathematics and}\\
{\footnotesize  Systems Science, Chinese Academy of Sciences,
Beijing 100190, China} }
\date{}

\usepackage{amssymb,amsmath,amsfonts,amsthm,color}

\def\D{\mathbb{D}}
\def\E{\mathbb{E}}

\def\R{\mathbb{R}}
\def\L{\mathcal{L}}
\def\S{\mathcal{S}}
\def\d{\textup{d}}

\def\Hess{\textup{Hess}}
\def\Ent{\textup{Ent}}

\def\Cylin{\textup{Cylin}}

\def\da{\downarrow}
\def\ra{\rightarrow}
\def\<{\langle}
\def\>{\rangle}
\def\ee{\varepsilon}
\def\fin{\hfill$\square$}

\newtheorem{theorem}{Theorem}[section]
\newtheorem{lemma}[theorem]{Lemma}

\newtheorem{proposition}[theorem]{Proposition}

\begin{document}

\maketitle
\makeatletter 
\renewcommand\theequation{\thesection.\arabic{equation}}
\@addtoreset{equation}{section}
\makeatother 

\begin{abstract}
We consider the Fokker-Planck equation on the abstract Wiener space
associated to the Ornstein-Uhlenbeck operator. Using the
Weitzenb\"ock formula, we prove an explicit estimate on the time
derivative of the entropy of the solution to the Fokker-Planck
equation.
\end{abstract}



\section{Introduction}

The second law in thermodynamics asserts that the entropy of an
isolated physical system always increases; an example is the entropy
of the solution to the heat equation on a Riemannian manifold.
Recently motivated by the ground-breaking work of Perelman, there
have been intensive studies on the entropy (i.e. Perelman's
$W$-functional) of the solution to the heat equation. In
\cite{Ni04a, Ni04b}, the author proved the monotonicity of the
entropy and showed its connection with the geometry of the manifold.
In the recent work \cite{LimLuo}, we presented several estimates on
the time derivative of entropy of the solution to the heat equation
on a Riemannian maniflod, in terms of the lower bound on the Ricci
curvature and the first eigenvalue of the Laplacian operator. In the
present work we intend to generalize these results to the infinite
dimensional case.

Let $(W,H,\mu)$ be an abstract Wiener space (see the beginning of
Section 2 for its definition) and $\L$ the Ornstein-Uhlenbeck
operator on $W$. Consider the Fokker-Planck equation
  \begin{equation}\label{sect-1.1}
  \frac\partial{\partial t}u_t=\L u_t,\quad u|_{t=0}=u_0,
  \end{equation}
where $u_0\in L^{p}(W)$ for some $p>1$. This equation is understood
in the weak sense, see \eqref{FPE.1}.  As shown in Section 3,
equation \eqref{sect-1.1} has a unique solution which is given by
$u_t=P_t u_0$ with $P_t$ being the Ornstein-Uhlenbeck semigroup on
$W$. Moreover $u_t$ also solves \eqref{sect-1.1} in the Fr\'echet
sense. Now suppose that $u_0>0$; then $u_t>0$ and $\int_W
u_t\,\d\mu\equiv \int_W u_0\,\d\mu$ for all $t\geq0$. Define the
entropy by
  $$\Ent(u_t)=-\int_W u_t\log u_t\,\d\mu.$$
From the simple inequality $-x\log x\leq 1-x$ for all $x\geq0$, we
know that $\Ent(u_t)\leq 1-\int_W u_0\,\d\mu$ for any $t>0$. The
formal calculation gives us
  \begin{align*}
  \frac{\d}{\d t}\Ent(u_t)&=-\int_W\bigg((\log u_t)\frac\partial{\partial t}u_t
  +\frac\partial{\partial t}u_t\bigg)\d\mu\cr
  &=-\int_W(\log u_t)\L u_t\,\d\mu
  =\int_W\frac{|\nabla u_t|_H^2}{u_t}\,\d\mu,
  \end{align*}
where the last equality follows from the integration by parts
formula. Therefore the entropy $\Ent(u_t)$ is an increasing function
of $t$ if the initial condition $u_0$ is not a constant. As in
\cite{LimLuo}, we will estimate the rate of change of the entropy as
$t\ra\infty$, by making use of the Weitzenb\"ock formula on $W$.

Denote by $\D_1^{p}(W)$ the first order Sobolev space on the Wiener
space $W$. The main result of this paper is

\begin{theorem}\label{main-thm} Let $p>1$. Suppose $u_0\in \D_1^{4p}(W)$ such
that $u_0\geq \ee_0$ for some positive constant $\ee_0>0$. Then
  $$\frac{\d}{\d t}\Ent(u_t)\leq e^{-2t}\int_W\frac{|\nabla u_0|_H^2}{u_0}\,\d\mu.$$
\end{theorem}

This theorem will be proved in Section 3. The above estimate is
consistent with the result in \cite[Example 2.4]{LimLuo}. Indeed,
let $\L_n$ be the $n$-dimensional version of the Ornstein-Uhlenbeck
operator $\L$ (see  \eqref{sect-2.3} for its definition); then we
have $\L_n=2L$ where $L$ is defined in \cite[Example 2.4]{LimLuo}
with $k=1$. Slight modification of the arguments in \cite[Example
2.4]{LimLuo} will give us that the time derivative of the entropy of
the transition density is $e^{-2t}$. We would like to mention that,
by \cite[Theorem 1.5]{FangShaoSturm09} and the metric measure theory
(see e.g. \cite{LottVillani, Sturm1, Sturm2}), the lower bound for
the Ricci curvature of the Wiener space $W$ is 1, hence the above
theorem is also in accordance with the main result in \cite{LimLuo}
if we consider the Laplacian operator $\Delta$ instead of
$\frac12\Delta$.

Compared to \cite[Theorem 1.1]{LimLuo}, the main difficulties in the
infinite dimensional situation are: (1) the justification of that
various functionals belong to the domain of the Ornstein-Uhlenbeck
operator $\L$, and (2) the validness of the differentiation under
the integral sign. Fortunately the unique solution $u_t$ to the
Fokker-Planck equation \eqref{sect-1.1} is sufficiently regular, and
the equation can actually be understood in the sense of Fr\'echet
differential. These observations make our computations possible. To
avoid the technical difficulties, we assume that the initial value
$u_0$ has a positive lower bound $\ee_0>0$, so that the estimations
become easier.

The paper is organized as follows. We recall in Section 2 some
preliminary elements in the Malliavin calculus and prove the
Weitzenb\"ock formula on $W$ associated to the Ornstein-Uhlenbeck
operator $\L$ (cf. Theorem \ref{sect-2-thm-1}). In Section 3 we
first show that $u_t=P_t u_0$ is the unique solution to the equation
\eqref{sect-1.1}; after that we establish an equality which is
essential for proving the main result of this paper, see Theorem
\ref{sect-3-thm-2}. Then by following the idea in \cite{LimLuo}, we
present the proof of Theorem \ref{main-thm}. Finally in the
Appendix, we give the proof of a result concerning the
differentiation under the integral sign which is needed in the proof
of the main theorem.

\section{Preliminaries in Malliavin calculus and the Weitzenb\"ock formula}

In this section we recall some basic facts in the Malliavin calculus
and present a Weitzenb\"ock type formula on the Wiener space
associated to the Ornstein-Uhlenbeck operator. Let $(W,H,\mu)$ be an
abstract Wiener space in the sense of L. Gross, i.e. $W$ is a
separable Banach space, $H$ is a separable Hilbert space and $\mu$
is a Borel probability on $W$, such that $H$ is continuously and
densely embedded into $W$ and for any $\ell\in W^\ast\,(\mbox{the
dual space of }W)$, we have
  $$\int_We^{\sqrt{-1}\,\ell(w)}\,\d\mu(w)=e^{-|\ell|_H^2/2},$$
where $|\cdot|_H$ is the norm in $H$ associated to the inner product
$\<\, ,\>_H$. In the following, we fix an orthonormal basis
$\{h_i:i\geq1\}$ of $H$, with $h_i\in W^\ast$ for all $i\geq1$.
Define $H_n=\mbox{span}\{h_i:1\leq i\leq n\}$ and
  $$\pi_n(w)=\sum_{i=1}^n h_i(w)\, h_i,\quad w\in W.$$
Then $\pi_n$ is a continuous linear map from $W$ onto $H_n$, and the
restriction of $\pi_n|_{H}$ is the orthogonal projection. It is
known that the push forward $\mu_n:=(\pi_n)_\#\mu$ is the standard
Gaussian measure on the $n$-dimensional Euclidean space $H_n$.

We refer to \cite{Fang05, IkedaWatanabe89, Malliavin97, Nualart} for
the background in Malliavin calculus. Let $K$ be a separable Hilbert
space. Denote by $K\otimes H$ the Hilbert space of Hilbert-Schmidt
operators $L$ from $K$ to $H$, and $\|L\|_{K\otimes H}$ the
Hilbert-Schmidt norm. For some $p>1$ and $Z\in L^p(W, K)$, we say
that $Z\in \D_1^p(W,K)$ if there exists $\nabla Z\in L^p(W, H\otimes
K)$ such that for each $h\in H$,
  $$\<\nabla Z,h\>_H=D_h Z=\frac{\d}{\d\ee}\Big|_{\ee=0}Z(w+\ee h)
  \quad \hbox{holds in }L^{p-}.$$
The space $\D_1^p(W,K)$ is complete under the norm:
  $$\|Z\|_{\D_1^p(W,K)}=\big(\|Z\|_{L^p(W,K)}^p + \|\nabla
  Z\|_{L^p(W,H\otimes K)}^p\big)^{1/p}.$$
In the same way we can define Sobolev spaces $\D_m^p(W,K)$ of higher
orders $m\geq1$. Notice that $\D_0^p(W,K)=L^p(W,K)$. A $K$-valued
functional $Z$ is called cylindrical if there exist $N,\,M\geq1$,
$f_i\in C^\infty_b(\R^M)$ and $k_i\in K\,(1\leq i\leq N)$, such that
  $$Z=\sum_{i=1}^N f_i(h_1(w),\cdots,h_M(w))\,k_i.$$
By the Schmidt orthogonalization procedure, we may always assume
that $\{k_1,\cdots,k_N\}$ is an orthonormal family. Note that
$Z:W\ra K$ is Fr\'{e}chet differentiable of any order. We denote by
$\Cylin(W,K)$ the space of $K$-valued cylindrical functionals, which
is dense in $\D_m^p(W,K)$. If $K=\R$, we simply write $\D^p_m(W)$
and $\Cylin(W)$. A basic result in Malliavin calculus is that the
divergence $\delta(Z)\in \D_{m-1}^p(W,K)$ exists for
$Z\in\D_m^p(W,H\otimes K)$ (see \cite[Proposition 1.5.7]{Nualart}),
and there is $C_{p,m}>0$ such that
  \begin{equation}\label{sect-2.1}
  \|\delta(Z)\|_{\D_{m-1}^p(W,K)}\leq C_{p,m}\, \|Z\|_{\D_m^p(W,H\otimes K)}.
  \end{equation}

The Mehler formula below defines the Ornstein-Uhlenbeck semigroup
$P_t$ on $W$:
  \begin{equation}\label{sect-2.2}
  P_t F(x)=\int_W
  F\big(e^{-t}x+\sqrt{1-e^{-2t}}\, y\big)\,\d\mu(y).
  \end{equation}
Here are some basic properties of $P_t$ that will be used later.

\begin{proposition}\label{OU-semigroup-property}
\begin{itemize}
\item[\rm(1)] For any $t>0$, $P_t(\Cylin(W))\subset \Cylin(W)$.

\item[\rm(2)] For any $t>0$ and $p\in[1,+\infty]$, we have for all $u\in
L^p(W)$, $\|P_t u\|_{L^p} \leq \|u\|_{L^p}$ and $\lim_{t\ra0} \|P_t
u-u\|_{L^p}=0$.

\item[\rm(3)] $P_\ee$ is self-adjoint in $L^2(W)$. Furthermore, for any
$p\in(1,+\infty)$ and $u\in L^p(W),\,v\in L^{q}(W)$ with
$\frac1p+\frac1{q}=1$, we have
  $$\int_W uP_\ee v\,\d\mu=\int_W vP_\ee u\,\d\mu.$$

\item[\rm(4)] For every $t>0$, $p>1$ and $m\geq1$, we have $P_t u\in
\D^p_m(W)$ for any $u\in L^p(W)$, and there is $C_{p,m}>0$ such that
  $$\|\nabla^m P_t u\|_{L^p(W,H^{\otimes m})}\leq C_{p,m}A_t^m\|u\|_{L^p(W)},$$
where $H^{\otimes m}=\underbrace{H\otimes\cdots\otimes H}_{m\
times}$ and $A_t=e^{-t}/\sqrt{1-e^{-2t}}$.
\end{itemize}
\end{proposition}

Recall that the last result is due to Sugita \cite{Sugita} (see also
\cite[Exercise 6.8]{Fang05}). The infinitesimal generator $\L$ of
$P_t$ is called the Ornstein-Uhlenbeck operator:
  \begin{equation}\label{sect-2.2.5}
  \frac{\partial}{\partial t}P_t F=\L P_t F=P_t \L F,\quad F\in\Cylin(W).
  \end{equation}
It is known that $\L=-\delta\circ\nabla$, hence $\L$ is a continuous
operator from $\D_{m+2}^p(W)$ to $\D_m^p(W)$ for all $m\geq 0$ and
$p>1$. Furthermore, if $F\in \Cylin(W)$ has the expression
$F(w)=f\circ\pi_n(w)$ for some $f\in C^\infty_b(H_n)$, then
  \begin{equation}\label{sect-2.3}
  \L F(w)=(\L_n f)\circ\pi_n(w),
  \end{equation}
where $\L_n$ is the Ornstein-Uhlenbeck operator on $H_n$: $\L_n
f(x)=\sum_{i=1}^n\big(\partial_{ii}f(x)-x^i\partial_i f(x)\big)$.
Our main result of this section is

\begin{theorem}[Weitzenb\"ock formula]\label{sect-2-thm-1}
Let $F\in \D_3^{2p}(W)$ for some $p>1$. Then
  \begin{equation}\label{sect-2-thm-1.1}
  \L(|\nabla F|_H^2)=2\<\nabla F,\nabla \L F\>_H+2|\nabla F|_H^2
  +2\|\nabla^2 F\|_{H\otimes H}^2.
  \end{equation}
\end{theorem}

An integral form of the Weitzenb\"ock formula is given in
\cite[Section 5.2]{Fang05}, which is the interpretation for 1-forms
of the de Rham-Hodge-Kodaira decomposition in \cite{Shigekawa}. The
equality \eqref{sect-2-thm-1.1} can be proved by following some of
the arguments in \cite[Section 5.3]{Fang05}. For the readers'
convenience, we include its proof here. First we consider the case
where $F$ is cylindrical.

\begin{lemma}\label{sect-2-lem-1}
The equality \eqref{sect-2-thm-1.1} holds for all $F\in \Cylin(W)$.
\end{lemma}

\noindent{\bf Proof.} Let $F(w)= f\circ\pi_n(w)$ for some smooth
function $f:H_n\ra\R$. We have $\nabla F(w)=\sum_{i=1}^n(\partial_i
f)(\pi_n(w))\, h_i$, hence
  \begin{equation}\label{sect-2-lem-1.1}
  |\nabla F(w)|_H^2=\sum_{i=1}^n\big[(\partial_i f)(\pi_n(w))\big]^2
  =|\nabla_n f|^2(\pi_n(w)),
  \end{equation}
where $\nabla_n$ is the gradient on the Euclidean space $H_n$. By
\eqref{sect-2.3}, we obtain
  \begin{equation}\label{sect-2-lem-1.2}
  \L\big(|\nabla F|_H^2\big)(w)=\big[\L_n\big(|\nabla_n f|^2\big)\big](\pi_n(w)).
  \end{equation}
Now direct computations lead to
  $$\L_n\big(|\nabla_n f|^2\big)(x)=2\|\Hess f\|_{H_n\otimes H_n}^2+2\<\nabla_n f,\nabla_n\Delta_n f\>
  -2(\Hess f)(x,\nabla_n f),\quad x\in \R^n,$$
where $\Delta_n$ is the Laplacian on $H_n$. Notice that $\nabla_n\<
x,\nabla_n f\>=\nabla_n f+(\Hess f)\cdot x$, therefore
  $$(\Hess f)(x,\nabla_n f)=\<\nabla_n\< x,\nabla_n f\>,\nabla_n f\>-|\nabla_n f|^2.$$
Substituting this into the above equality and by the definition of
$\L_n$, we get
  \begin{align*}
  \L_n\big(|\nabla_n f|^2\big)(x)&=2\|\Hess f\|_{H_n\otimes H_n}^2+2\<\nabla_n f,\nabla_n\Delta_n f\>
  -2\<\nabla_n\< x,\nabla_n f\>,\nabla_n f\>+2|\nabla_n f|^2\cr
  &=2\|\Hess f\|_{H_n\otimes H_n}^2+2\<\nabla_n f,\nabla_n\L_n f\>+2|\nabla_n
  f|^2.
  \end{align*}
Combining this with \eqref{sect-2.3}, \eqref{sect-2-lem-1.1} and
\eqref{sect-2-lem-1.2}, we get the desired result. \fin

\medskip

\noindent{\bf Proof of Theorem \ref{sect-2-thm-1}.} Since
$\Cylin(W)$ is dense in $\D_3^{2p}(W)$, there exists a sequence
$\{F_n:n\geq1\}$ of cylindrical functionals such that
$\lim_{n\ra\infty} \|F_n-F\|_{\D_3^{2p}(W)}=0$. By Lemma
\ref{sect-2-lem-1}, for all $n\geq1$,
  \begin{equation}\label{sect-2.4}
  \L(|\nabla F_n|_H^2)=2\<\nabla F_n,\nabla \L F_n\>_H+2|\nabla F_n|_H^2
  +2\|\nabla^2 F_n\|_{H\otimes H}^2.
  \end{equation}
It remains to show that both sides of the above equality converge in
$L^p(W)$. First it is clear that $|\nabla F_n|_H^2$ converges to
$|\nabla F|_H^2$ in $\D_2^p(W)$. Note that the Ornstein-Uhlenbeck
operator is a continuous map from $\D_2^p(W)$ to $L^p(W)$, thus
  \begin{equation}\label{sect-2.5}
  \lim_{n\ra\infty}\big\|\L(|\nabla F_n|_H^2)-\L(|\nabla
  F|_H^2)\big\|_{L^p}=0.
  \end{equation}
Next by the triangular inequality,
  \begin{align*}
  \|\<\nabla F_n,\nabla \L F_n\>_H-\<\nabla F,\nabla \L F\>_H\|_{L^p}
  &\leq \|\<\nabla F_n-\nabla F,\nabla \L F_n\>_H\|_{L^p}\cr
  &\hskip12pt +\|\<\nabla F,\nabla \L F_n-\nabla \L F\>_H\|_{L^p}.
  \end{align*}
Cauchy's inequality leads to
  \begin{align}\label{sect-2.6}
  \|\<\nabla F_n-\nabla F,\nabla \L F_n\>_H\|_{L^p}& \leq \|\nabla F_n-\nabla F\|_{L^{2p}(W,H)}
  \|\nabla \L F_n\|_{L^{2p}(W,H)}\cr
  &\leq \|F_n-F\|_{\D_1^{2p}(W)}\cdot C_p\|F_n\|_{\D_3^{2p}(W)},
  \end{align}
where the last inequality follows from the boundedness of the
operator $\L:\D_3^{2p}(W)\ra \D_1^{2p}(W)$. Since the sequence
$\|F_n\|_{\D_3^{2p}(W)}$ is bounded, we arrive at
  $$\lim_{n\ra\infty}\|\<\nabla F_n-\nabla F,\nabla \L F_n\>_H\|_{L^p}=0.$$
Similarly we have
  $$\|\<\nabla F,\nabla \L F_n-\nabla \L F\>_H\|_{L^p}\leq
  C_p\|F\|_{\D_1^{2p}(W)}\|F_n-F\|_{\D_3^{2p}(W)}$$
whose right hand side tends to 0 as $n\ra\infty$. From this and
\eqref{sect-2.6}, we conclude that the first term on the right hand
side of \eqref{sect-2.4} converges in $L^p(W)$ to $2\<\nabla
F,\nabla \L F\>_H$. Finally it is easy to show that the last two
terms of equality \eqref{sect-2.4} also converge in $L^p(W)$ to
$2|\nabla F|_H^2$ and $2\|\nabla^2 F\|_{H\otimes H}^2$,
respectively. Combining these results with \eqref{sect-2.5}, we
complete the proof of Theorem \ref{sect-2-thm-1} by taking limit in
\eqref{sect-2.4}. \fin

\section{Entropy of the solution to the Fokker-Planck equation on $W$}

In this section we will estimate the time derivative of entropy of
the solution to the Fokker-Planck equation associated to the
Ornstein-Uhlenbeck operator. From \eqref{sect-2.2.5}, we see that
for any initial value $u_0\in \Cylin(W)$, $u_t:=P_t u_0$ gives the
solution to the Fokker-Planck equation in the classical sense. In
the following we show that this fact remains true for all $u_0\in
L^p(W)$, where $p>1$.

First we introduce the notion of weak solution to \eqref{sect-1.1}.
A function $u\in L^\infty([0,\infty),L^p(W))$ is called a weak
solution to the equation \eqref{sect-1.1} if for any $\alpha\in
C_c^\infty([0,\infty))$ and $F\in \Cylin(W)$, it holds
  \begin{equation}\label{FPE.1}
  -\alpha(0)\int_W Fu_0\,\d\mu=\int_0^\infty\!\!\int_W\big[\alpha^\prime(t)F+\alpha(t)\L
  F\big]u_t\,\d\mu\d t.
  \end{equation}
For $u_0\geq 0$ satisfying $\int_W u_0\,\d\mu=1$ and
$-\Ent(u_0)=\int_W u_0\log u_0\,\d\mu<+\infty$, it is shown in
\cite[Theorem 3.7]{FangShaoSturm09} that the weak solution $u_t$ of
equation \eqref{sect-1.1} can be constructed via De Giorgi's
``minimizing movement'' approximation. But the uniqueness of
solutions is not considered there. The existence and uniqueness of
general Fokker-Planck type equations are studied in \cite{Luo10};
however, those results do not apply to the equation \eqref{FPE.1}
(see \cite[Remark 4.6]{Luo10}).

We first prove the following simple result. For $n\geq1$, let
$\E^{H_n}$ be the conditional expectation on $W$ with respect to the
$\sigma$-field generated by cylindrical functionals of the form
$F=f\circ\pi_n$. We also denote by $\E^\mu$ and $\E^{\mu_n}$ the
expectations on $W$ and $H_n$ respectively.

\begin{lemma}\label{sect-3-lem-1}
Let $u_0\in L^p(W)$ and $u_t=P_tu_0,\, t\geq0$. For $n\geq1$, define
the function $u_n(t)\in L^p(H_n,\mu_n)$ such that
$u_n(t)\circ\pi_n=\E^{H_n}(u_t)$. Then
  $$u_n(t)=P^{(n)}_t u_n(0),$$
where $P^{(n)}_t$ is the Ornstein-Uhlenbeck semigroup on $H_n$.
\end{lemma}

\noindent{\bf Proof.} For any $F=f\circ\pi_n$ with $f\in
C_b^\infty(H_n)$, by the symmetry of the operator $P_t$ (see
Proposition \ref{OU-semigroup-property}(3)), we have
  \begin{align*}
  \E^{\mu_n}[fu_n(t)]&=\E^\mu\big[(f\circ\pi_n)\,\E^{H_n}(u_t)\big]
  =\E^\mu[F u_t]=\E^\mu[(P_t F)u_0].
  \end{align*}
By the Mehler formula \eqref{sect-2.2}, it is clear that $P_t
F=\big(P^{(n)}_t f\big)\circ\pi_n$, hence
  \begin{align*}
  \E^{\mu_n}[fu_n(t)]&=\E^\mu\big[\big(P^{(n)}_t f\big)\circ\pi_n\cdot u_0\big]
  =\E^\mu\big[\big(P^{(n)}_t f\big)\circ\pi_n\cdot
  \E^{H_n}(u_0)\big]\cr
  &=\E^{\mu_n}\big[\big(P^{(n)}_t f\big)\cdot u_n(0)\big]
  =\E^{\mu_n}\big[f P^{(n)}_t u_n(0)\big],
  \end{align*}
where the last equality follows from the symmetry of $P^{(n)}_t$.
Since $f\in C_b^\infty(H_n)$ is arbitrary, we complete the proof.
\fin

\begin{theorem}\label{sect-3-thm-1}
Let $u_0\in L^p(W)$ with $p>1$. Then $u_t:=P_t u_0$ is the unique
weak solution to the Fokker-Planck equation \eqref{sect-1.1} in the
space
  $$\S=\big\{v\in L^\infty\big([0,\infty),L^p(W)\big):
  \nabla v\in L^1\big([0,\infty),L^p(W,H)\big)\big\}.$$
\end{theorem}

\noindent{\bf Proof.} First we check that $u_t=P_t u_0$ is really a
solution of \eqref{sect-1.1}. Letting $m=1$ in Proposition
\ref{OU-semigroup-property}(4) and by the definition of $A_t$, it is
clear that $\nabla(P_tu_0)\in L^1\big([0,\infty),L^p(W,H)\big)$,
hence $u\in\S$. Next we take a sequence $\{u^{(n)}:n\geq1\}$ of
cylindrical functionals such that $\lim_{n\ra\infty}
\|u^{(n)}-u_0\|_{L^p}=0$. For any $n\geq1$, let $u^{(n)}_t=P_t
u^{(n)}$. Then by \eqref{sect-2.2.5}, we have
$\frac{\partial}{\partial t}u^{(n)}_t= \L u^{(n)}_t$. Thus for any
$\alpha\in C_c^\infty([0,\infty))$ and $F\in \Cylin(W)$, it holds
  \begin{equation}\label{sect-3-thm-1.1}
  -\alpha(0)\int_W Fu^{(n)}\,\d\mu=\int_0^\infty\!\!\int_W\big[\alpha^\prime(t)F+\alpha(t)\L
  F\big]u^{(n)}_t\,\d\mu\d t.
  \end{equation}
By Proposition \ref{OU-semigroup-property}(2), we have
  $$\sup_{0\leq t<\infty}\|u^{(n)}_t-u_t\|_{L^p}\leq
  \|u^{(n)}-u_0\|_{L^p}\ra0 $$
as $n\ra\infty$. Noting that $\L F\in L^q(W)$ where $q$ is the
conjugate number of $p:p^{-1}+q^{-1}=1$, taking limit in
\eqref{sect-3-thm-1.1} gives us the equality \eqref{FPE.1}. That is
to say, $u_t=P_t u_0$ is a weak solution to the Fokker-Planck
equation \eqref{sect-1.1}.

Next we prove the uniqueness of weak solutions to \eqref{FPE.1}. Let
$v\in \S$ be any solution of \eqref{FPE.1} with $v_0=u_0$. For any
$n\geq1$, define $u_n\in L^p(H_n,\mu_n)$ such that
  $$u_n\circ\pi_n=\E^{H_n}(u_0).$$
Similarly we define $v_n(t)\in L^p(H_n,\mu_n)$ from $v_t$. It is
easy to know that $v_n$ belongs to the space
  $$\S_n=\big\{v_n\in L^\infty\big([0,\infty),L^p(H_n,\mu_n)\big):
  \nabla_n v_n\in L^1\big([0,\infty),L^p(H_n,H_n,\mu_n)\big)\big\}.$$
Now fix $n\geq 1$. Then for any $F=f\circ\pi_n\in \Cylin(W)$ and
$\alpha\in C_c^\infty([0,\infty))$, it holds
  \begin{equation}\label{sect-3-thm-1.2}
  -\alpha(0)\int_W Fu_0\,\d\mu=\int_0^\infty\!\!\int_W[\alpha^\prime(t)F+\alpha(t)\L
  F]v_t\,\d\mu\d t.
  \end{equation}
We have
  \begin{align*}
  \int_W Fu_0\,\d\mu&=\E^\mu\big[(f\circ\pi_n)u_0\big]
  =\E^\mu\big[(f\circ\pi_n)\E^{H_n}(u_0)\big]\cr
  &=\E^{\mu}\big[(f\circ\pi_n)(u_n\circ\pi_n)\big]
  =\int_{H_n}fu_n\,\d\mu_n.
  \end{align*}
Now for any $t\in[0,T]$, in the same way we have
  \begin{align*}
  \int_W Fv_t\,\d\mu=\int_{H_n}fv_n(t)\,\d\mu_n
  \end{align*}
and by \eqref{sect-2.3},
  $$ \int_W (\L F)v_t\,\d\mu=\int_{H_n}(\L_n f) v_n(t)\,\d\mu_n.$$
Therefore the equation \eqref{sect-3-thm-1.2} becomes
  \begin{equation}\label{sect-3-thm-1.3}
  -\alpha(0)\int_{H_n}fu_n\,\d\mu_n=\int_0^\infty\!\!\int_{H_n}\big[\alpha^\prime(t)f+\alpha(t)\L_n
  f\big]v_n(t)\,\d\mu_n\d t.
  \end{equation}
This means that $v_n$ is a weak solution to the finite dimensional
Fokker-Planck equation
  $$\frac{\partial}{\partial t}v_n(t)=\L_n v_n(t)$$
with initial condition $v_n(0)=u_n$. Approximating $u_n$ by smooth
functions with respect to the norm in $L^p(H_n,\mu_n)$ and following
the argument of the starting part of this theorem, we can show that
$P^{(n)}_t u_n$ is also the solution to \eqref{sect-3-thm-1.3} with
the same initial value $u_n$. Moreover $P^{(n)}_t u_n\in \S_n$. By
the uniqueness of solutions in the finite dimensional case (cf.
\cite[Theorem 4.10]{Luo10} restricted to the finite dimensional
context or \cite[Corollary 1]{LeBrisLions08}), we obtain
$v_n(t)=P^{(n)}_t u_n$ for all $t>0$. Lemma \ref{sect-3-lem-1} tells
us that $v_n(t)\circ\pi_n=\big(P^{(n)}_t u_n\big)\circ\pi_n=
\E^{H_n}(u_t)$. Letting $n$ tend to infinity, we conclude that
$v_t=u_t=P_t u_0$. Therefore we get the uniqueness of weak solutions
to the Fokker-Planck equation \eqref{FPE.1}. \fin

\medskip

In the next proposition we show that $u_t=P_t u_0$ is also a
solution to equation \eqref{sect-1.1} in the strong sense.

\begin{proposition}\label{sect-3-prop-1}
Let $u_0\in L^p(W)$ for some $p>1$. Then for any $m\geq 0$, the
following equality holds in $\D_m^p(W)$:
  $$\frac{\partial}{\partial t}P_t u_0=\L(P_tu_0),\quad \mbox{for all }t>0.$$
\end{proposition}

\noindent{\bf Proof.} By Proposition \ref{OU-semigroup-property}(4),
for any $t>0$,
  \begin{align}\label{sect-3-prop-1.1}
  \|P_t u_0\|_{\D_m^p(W)}\leq \bigg[\sum_{i=0}^m
  C_{p,i}^pA_t^{ip}\|u_0\|_{L^p}^p\bigg]^{1/p}
  \leq \bar C_{p,m}\|u_0\|_{L^p}\sum_{i=0}^m A_t^i.
  \end{align}
Therefore by the boundedness of the Ornstein-Uhlenbeck operator
$\L:\D_{m+2}^p(W) \ra \D_m^p(W)$,
  \begin{align}\label{sect-3-prop-1.2}
  \|\L P_t u_0\|_{\D_m^p(W)}\leq \hat C_{p,m}\|P_t u_0\|_{\D_{m+2}^p(W)}
  \leq \tilde C_{p,m}\|u_0\|_{L^p}\sum_{i=0}^{m+2} A_t^i.
  \end{align}
Recall that $A_t=e^{-t}/\sqrt{1-e^{-2t}}$. Thus for any $0<s<t<
\infty$, the right hand side is integrable on the interval $[s,t]$.
In particular, \eqref{sect-3-prop-1.1} implies that the curve
$(0,\infty)\ni t\ra P_t u_0$ is locally integrable in the Sobolev
space $\D_m^p(W)$.

Now we take a sequence $\{u_n:n\geq1\}$ of cylindrical functionals
such that $\|u_n-u_0\|_{L^p(W)}\ra0$ as $n$ goes to $\infty$. Then
for any $t>0$, \eqref{sect-3-prop-1.1} leads to
  \begin{align}\label{sect-3-prop-1.3}
  \varlimsup_{n\ra\infty}\|P_t u_n-P_t u_0\|_{\D_m^p(W)}\leq
  \bar C_{p,m}\bigg(\sum_{i=0}^m A_t^i\bigg)\lim_{n\ra\infty}\|u_n-u_0\|_{L^p}=0.
  \end{align}
Similarly, for all $0<s<t<\infty$, by \eqref{sect-3-prop-1.2}
  \begin{align}\label{sect-3-prop-1.4}
  \bigg\|\int_s^t\L P_\tau u_n\,\d\tau-\int_s^t\L P_\tau u_0\,\d\tau\bigg\|_{\D_m^p(W)}
  &\leq \int_s^t\|\L P_\tau u_n-\L P_\tau u_0\|_{\D_m^p(W)}\,\d\tau\cr
  &\leq \tilde C_{p,m}\|u_n-u_0\|_{L^p}\sum_{i=0}^{m+2}\int_s^t
  A_\tau^i\,\d\tau,
  \end{align}
whose right hand side tends to 0 as $n\ra\infty$. Now by
\eqref{sect-2.2.5}, for every $n\geq1$, we have
  $$P_t u_n-P_s u_n=\int_s^t \L P_\tau u_n\,\d\tau.$$
With \eqref{sect-3-prop-1.3} and \eqref{sect-3-prop-1.4} in mind,
letting $n\ra\infty$ in the above equality gives us
  \begin{equation}\label{sect-3-prop-1.5}
  P_t u_0-P_s u_0=\int_s^t \L P_\tau u_0\,\d\tau,
  \end{equation}
which holds in any Sobolev space $\D_m^p(W)$. In particular, by
\eqref{sect-3-prop-1.2},
  $$\|P_t u_0-P_s u_0\|_{\D_m^p(W)}\leq \int_s^t \|\L P_\tau u_0\|_{\D_m^p(W)}\,\d\tau
  \leq \tilde C_{p,m}\|u_0\|_{L^p}\sum_{i=0}^{m+2}\int_s^t A_\tau^i\,\d\tau,$$
which implies that $(0,\infty)\ni t\mapsto P_t u_0\in\D_m^p(W)$ is
continuous for any $m\geq1$.

Now by the boundedness of $\L$,
  \begin{align*}
  \bigg\|\frac{P_t u_0-P_s u_0}{t-s}-\L P_s u_0\bigg\|_{\D_m^p(W)}
  &\leq \frac 1{t-s}\int_s^t\|\L P_\tau u_0-\L P_s
  u_0\|_{\D_m^p(W)}\,\d\tau\cr
  &\leq \frac {C^\prime_{p,m}}{t-s}\int_s^t\|P_\tau u_0-P_s
  u_0\|_{\D_{m+2}^p(W)}\,\d\tau.
  \end{align*}
The continuity of $(0,\infty)\ni \tau\mapsto P_\tau
u_0\in\D_{m+2}^p(W)$ gives rise to
  $$\lim_{t\ra s}\bigg\|\frac{P_t u_0-P_s u_0}{t-s}-\L P_s u_0\bigg\|_{\D_m^p(W)}=0.$$
The proof is complete. \fin

\medskip

Now we prove an equality which is critical in the proof of the main
result.

\begin{theorem}\label{sect-3-thm-2}
Let $u_t=P_t u_0$ with $u_0\in L^{4p}(W)$ and $u_0\geq \ee_0$ for
some $\ee_0>0$. Then
  $$\bigg(\L-\frac{\partial}{\partial t}\bigg)\bigg(\frac{|\nabla u_t|_H^2}{u_t}\bigg)
  =\frac2{u_t}|\nabla u_t|_H^2+\frac 2{u_t}\bigg\|\nabla^2 u_t
  -\frac{\nabla u_t\otimes\nabla u_t}{u_t}\bigg\|_{H\otimes H}^2.$$
\end{theorem}

\noindent{\bf Proof.} By Proposition \ref{OU-semigroup-property}(4),
$u_t=P_t u_0\in \D_3^{4p}(W)$, thus $|\nabla
u_t|_H^2\in\D_2^{2p}(W)$. We also have $u_t^{-1}\in \D_2^{2p}(W)$.
Indeed, since the initial value $u_0$ is bounded from below by
$\ee_0>0$, we have $u_t\geq\ee_0$ for all $t\geq0$. From
$\nabla(u_t^{-1})=-u_t^{-2}\nabla u_t$ it follows that
$|\nabla(u_t^{-1})|_H\leq \ee_0^{-2}|\nabla u_t|_H$. Thus
$\nabla(u_t^{-1}) \in L^{4p}(W,H)$. Next
  $$\nabla^2(u_t^{-1})=2u_t^{-3}\nabla u_t\otimes \nabla u_t-u_t^{-2}\nabla^2 u_t,$$
hence
  $$\|\nabla^2(u_t^{-1})\|_{H\otimes H}\leq 2\ee_0^{-3}|\nabla u_t|_H^2
  +\ee_0^{-2}\|\nabla^2 u_t\|_{H\otimes H}$$
which implies that $\nabla^2(u_t^{-1})\in L^{2p}(W,H\otimes H)$. To
sum up, $u_t^{-1}\in \D_2^{2p}(W)$. Therefore by \cite[Proposition
1.5.6]{Nualart}, we see that $\frac{|\nabla u_t|_H^2}{u_t}\in
\D_2^p(W)$.

Now
  \begin{equation}\label{sect-3-thm-2.1}
  \L\bigg(\frac{|\nabla u_t|^2}{u_t}\bigg)=u_t^{-1}\L(|\nabla u_t|^2)
  +|\nabla u_t|^2\L(u_t^{-1})+2\big\<\nabla(u_t^{-1}),\nabla(|\nabla
  u_t|^2)\big\>_H.
  \end{equation}
By the Weitzenb\"ock formula proved in Theorem \ref{sect-2-thm-1},
  \begin{equation}\label{sect-3-thm-2.2}
  \L(|\nabla u_t|_H^2)=2\<\nabla u_t,\nabla \L u_t\>_H+2|\nabla u_t|_H^2
  +2\|\nabla^2 u_t\|_{H\otimes H}^2.
  \end{equation}
It is easy to show that $\L(u_t^{-1})=-u_t^{-2}\L u_t+2u_t^{-3}
|\nabla u_t|_H^2$ and
  $$\big\<\nabla(u_t^{-1}),\nabla(|\nabla u_t|^2)\big\>_H
  =-2u_t^{-2}\big\<\nabla^2u_t,\nabla u_t\otimes\nabla u_t\big\>_{H\otimes H}.$$
Substituting these equalities and \eqref{sect-3-thm-2.2} into
\eqref{sect-3-thm-2.1} gives us
  $$\L\bigg(\frac{|\nabla u_t|_H^2}{u_t}\bigg)=\frac2{u_t}\big\<\nabla u_t,\nabla \L u_t\big\>_H
  -\frac{|\nabla u_t|_H^2}{u_t^2}\L u_t+\frac2{u_t} |\nabla u_t|_H^2
  +\frac 2{u_t}\bigg\|\nabla^2 u_t
  -\frac{\nabla u_t\otimes\nabla u_t}{u_t}\bigg\|_{H\otimes H}^2.$$
Next it is clear that
  \begin{equation}\label{sect-3-thm-2.3}
  \frac{\partial}{\partial t}\bigg(\frac{|\nabla u_t|_H^2}{u_t}\bigg)
  =\frac2{u_t}\big\<\nabla u_t,\nabla \big(\frac{\partial}{\partial t} u_t\big)\big\>_H
  -\frac{|\nabla u_t|_H^2}{u_t^2}\frac{\partial}{\partial t} u_t.
  \end{equation}
Combining the above two equalities, we obtain the desired formula.
\fin

\medskip

Finally we are in the position to prove the main result of this
paper.

\medskip

\noindent{\bf Proof of Theorem \ref{main-thm}.} By the definition of
entropy,
  $$\frac{\d}{\d t}\Ent(u_t)=-\frac{\d}{\d t}\int_Wu_t\log u_t\,\d\mu.$$
To commute the differential and integral, we have to check the
conditions in Theorem \ref{sect-3-thm-3} for the function
$(0,\infty)\times W\ni(t,w)\mapsto u_t(w)\log u_t(w)$. The first
condition is obviously satisfied. In view of Proposition
\ref{sect-3-prop-1}, we have
  $$\frac{\partial}{\partial t}(u_t\log u_t)=(\log u_t)\L u_t+\L u_t,$$
hence the condition (ii) is verified. It remains to check condition
(iii). Let $0<a<b<\infty$. Taking $m=0$ in \eqref{sect-3-prop-1.2},
we obtain that
  \begin{equation}\label{Proof.1}
  \|\L u_t\|_{L^{2p}(W)}\leq \tilde C_p\|u_0\|_{L^{2p}}\sum_{i=0}^2
  A_t^i,\quad t>0.
  \end{equation}
Since $u_t\geq\ee_0>0$, it is clear that $\log\ee_0\leq \log u_t\leq
u_t$, hence
  $$|\log u_t|\leq |\log \ee_0|\vee u_t\leq |\log \ee_0|+u_t.$$
Therefore by the contraction property of the Ornstein-Uhlenbeck
semigroup $P_t$,
  \begin{equation}\label{Proof.2}
  \sup_{t>0}\|\log u_t\|_{L^{2p}}\leq |\log \ee_0|+\sup_{t>0}\|u_t\|_{L^{2p}}
  \leq |\log \ee_0|+\|u_0\|_{L^{2p}}.
  \end{equation}
Now by \eqref{Proof.1} and \eqref{Proof.2}, Cauchy's inequality
gives us
  $$\int_a^b\!\!\int_W\big|(\log u_t)\L u_t\big|\,\d\mu\d t
  \leq \int_a^b\|(\log u_t)\L u_t\|_{L^p}\,\d t
  \leq \int_a^b\|\log u_t\|_{L^{2p}}\|\L u_t\|_{L^{2p}}\,\d t<+\infty.$$
Hence the condition (iii) in Theorem \ref{sect-3-thm-3} is satisfied
too. We pass the differentiation into the integral sign and get
  $$\frac{\d}{\d t}\Ent(u_t)=-\int_W\big[(\log u_t)\L u_t+\L u_t\big]\d\mu.$$

We have
  $$|\nabla\log u_t|_H=|u_t^{-1}\nabla u_t|_H\leq \ee_0^{-1}|\nabla u_t|_H,$$
hence $\log u_t\in \D_1^{2p}(W)$. By the integration by parts
formula we get
  \begin{equation}\label{Proof.3}
  \frac{\d}{\d t}\Ent(u_t)=-\int_W(\log u_t)\L u_t\,\d\mu
  =\int_W\frac{|\nabla u_t|_H^2}{u_t}\,\d\mu.
  \end{equation}
Next we have $\int_W \L\big(\frac{|\nabla
u_t|^2}{u_t}\big)\,\d\mu=0$, again due to the integration by parts
formula. By \eqref{sect-3-thm-2.3}, we can prove in a similar way
that the function $(t,w)\mapsto \frac{|\nabla u_t|_H^2}{u_t} (w)$
satisfies the three conditions in Theorem \ref{sect-3-thm-3}, hence
  $$\frac{\d}{\d t}\int_W\frac{|\nabla u_t|_H^2}{u_t}\,\d\mu
  =\int_W\frac{\partial}{\partial t}\bigg(\frac{|\nabla u_t|_H^2}{u_t}\bigg)\d\mu.$$
Therefore, integrating both sides of the formula in Theorem
\ref{sect-3-thm-2}, we obtain
  \begin{align*}
  -\frac{\d}{\d t}\int_W\frac{|\nabla u_t|_H^2}{u_t}\,\d\mu
  &=2\int_W\frac1{u_t}\bigg|\nabla^2 u_t
  -\frac{\nabla u_t\otimes\nabla u_t}{u_t}\bigg|_{H\otimes H}^2\,\d\mu
  +2\int_W\frac{|\nabla u_t|_H^2}{u_t}\,\d\mu\cr
  &\geq2\int_W\frac{|\nabla u_t|_H^2}{u_t}\,\d\mu.
  \end{align*}
Using Proposition \ref{OU-semigroup-property}(2) (it is also true
for functionals in $L^p(W,H)$), we can show that $t\mapsto
\int_W\frac{|\nabla u_t|_H^2}{u_t}\,\d\mu$ is right continuous at
$t=0$. Therefore
  $$\int_W\frac{|\nabla u_t|_H^2}{u_t}\,\d\mu\leq e^{-2t}\int_W\frac{|\nabla u_0|_H^2}{u_0}\,\d\mu.$$
In view of \eqref{Proof.3}, the proof is complete. \fin

\section{Appendix: a result on the differentiation under the integral sign}

In the proof of Theorem \ref{main-thm}, we need the following
theorem from the analysis which guarantees the differentiation under
the integral sign. See \cite{Cheng} for an introduction of related
results. For the reader's convenience, we give its complete proof
here.

\begin{theorem}\label{sect-3-thm-3}
Let $T$ be an open interval of $\R$, and $(\Omega,\nu)$ a measure
space. Suppose that a function $f:T\times\Omega\ra\R$ satisfies the
following conditions:
\begin{itemize}
\item[\rm(i)] $f(t,\omega)$ is a measurable function of $t$ and $\omega$
jointly, and is integrable over $\Omega$ for almost all $t\in T$
fixed;
\item[\rm(ii)] for almost all $\omega$, the derivative $\frac{\partial}
{\partial t}f(t,\omega)$ exists for all $t \in T$;
\item[\rm(iii)] for all compact intervals $[a,b]\subset T$, we have
  $$\int_a^b\!\!\int_\Omega \bigg|\frac{\partial}
  {\partial t}f(t,\omega)\bigg|\d\nu(\omega)\d t<+\infty.$$
\end{itemize}
Then for a.e. $t\in T$,
  $$\frac{\d}{\d t}\int_\Omega f(t,\omega)\,\d\nu(\omega)
  =\int_\Omega \frac{\partial}
  {\partial t}f(t,\omega)\,\d\nu(\omega).$$
\end{theorem}

\noindent{\bf Proof.} For $t\in T$, define $F(t)=\int_\Omega
f(t,\omega)\,\d\nu(\omega)$. Then $t\mapsto F(t)$ is absolutely
continuous on $T$. Indeed, for any $\ee>0$, we conclude from
condition (iii) that there is $\kappa>0$, such that for any
measurable set $E\subset [a,b]$ whose Lebesgue measure is less than
$\kappa$, it holds
  $$\int_E\!\int_\Omega \bigg|\frac{\partial}
  {\partial t}f(t,\omega)\bigg|\d\nu(\omega)\d t<\ee.$$
For any finite sequence of pairwise disjoint intervals $(x_k, y_k)$
of $T$ satisfying $\cup_k(x_k, y_k)\subset [a,b]$ and
$\sum_k(y_k-x_k)<\kappa$, we have
  \begin{align}\label{sect-3-thm-3.0}
  \sum_k|F(y_k)-F(x_k)|\leq
  \sum_k\int_\Omega|f(y_k,\omega)-f(x_k,\omega)|\,\d\nu(\omega).
  \end{align}
By (iii) and Fubini's theorem, for any $0<a<b<\infty$,
  $$\int_\Omega\!\int_a^b\bigg|\frac{\partial}
  {\partial t}f(t,\omega)\bigg|\,\d s\d\nu(\omega)<\infty,$$
hence for $\nu$-a.e. $\omega\in\Omega$, $\int_a^b
\big|\frac{\partial}{\partial t}f(t,\omega)\big|\,\d t<\infty$. As a
result,
  $$f(b,\omega)-f(a,\omega)=\int_a^b\frac{\partial}{\partial t}f(t,\omega)\,\d t.$$
By \eqref{sect-3-thm-3.0} and Fubini's theorem, we obtain
  \begin{align*}
  \sum_k|F(y_k)-F(x_k)|&\leq
  \sum_k\int_\Omega\!\int_{x_k}^{y_k}\bigg|\frac{\partial}
  {\partial t}f(t,\omega)\bigg|\,\d s\d\nu(\omega)\cr
  &=\int_{\cup_k (x_k,y_k)}\!\int_\Omega\bigg|\frac{\partial}
  {\partial t}f(t,\omega)\bigg|\,\d\nu(\omega)\d s<\ee.
  \end{align*}
Therefore $F$ is locally absolutely continuous, and $\frac{\d}{\d
t}F(t)$ exists for a.e. $t\in T$.

Next for $t\in T$, define $g(t)=\int_\Omega \frac{\partial}
{\partial t}f(t,\omega)\,\d\nu(\omega)$. By assumption (iii), the
function $g$ is locally integrable on $T$. Fix $a\in T$, define
$G(t)=\int_a^t g(s)\,\d s$ for all $t\in T$. Then by Lebesgue's
differentiation theorem, we have
  \begin{equation}\label{sect-3-thm-3.1}
  g(t)=\frac{\d}{\d t}G(t)\quad \mbox{for a.e. }t\in T.
  \end{equation}
Next for a.e. $t\in T$,
  \begin{align*}
  \frac{\d}{\d t}G(t)=\lim_{h\da0}\frac1h\big(G(t+h)-G(t)\big)
  =\lim_{h\da0}\frac1h\int_t^{t+h}\!\!\int_\Omega\frac{\partial}{\partial
  s}f(s,\omega)\,\d\nu(\omega)\d s.
  \end{align*}
Due to condition (iii), we can apply Fubini's theorem to get
  \begin{align*}
  \frac{\d}{\d t}G(t)&=\lim_{h\da0}\frac1h\int_\Omega\!\int_t^{t+h}\frac{\partial}{\partial
  s}f(s,\omega)\,\d s\d\nu(\omega)\cr
  &=\lim_{h\da0}\frac1h\int_\Omega\big(f(t+h,\omega)-f(t,\omega)\big)\,\d\nu(\omega)
  =\frac{\d}{\d t}F(t).
  \end{align*}
Combining this with \eqref{sect-3-thm-3.1} and the definitions of
$F,g$, we complete the proof. \fin


\begin{thebibliography}{a23}

\bibitem{Cheng} S. Cheng, {\it Differentiation under the integral
sign}, http://planetmath.org/encyclopedia/Dif\newline
ferentiationUnderIntegralSign.html\#tex2html13.

\bibitem{Fang05} S. Fang, {\it Introduction to Malliavin calculus}.
Math. Series for Graduate students, vol. 3, 2005, Tsinghua
University Press, Springer.

\bibitem{FangShaoSturm09} S. Fang, J. Shao and K.T. Sturm, {\it Wasserstein
space over the Wiener space}. Probab. Theory Relat. Fields 146
(2010), no. 3-4, 535--565.

\bibitem{IkedaWatanabe89} N. Ikeda and S. Watanabe, {\it Stochastic Differential Equations
and Diffusion Processes}. North-Holland, Amsterdam, 1989.

\bibitem{LeBrisLions08} C. LeBris and P.L. Lions, {\it Existence and
uniqueness of solutions to Fokker-Planck type equations with
irregular coefficients}. Comm. Partial Differential Equations 33
(2008), 1272--1317.

\bibitem{LimLuo} Adrian P. C. Lim and Dejun Luo, {\it Asymptotic estimates
on the time derivative of entropy on a Riemannian manifold}.
arXiv:1011.3979.

\bibitem{LottVillani} J. Lott and C. Villani, {\it Ricci curvature for metric-measure
spaces via optimal transport}. Ann. of Math. 169 (2009), no. 3,
903--991.

\bibitem{Luo10} Dejun Luo, {\it Well-posedness of Fokker-Planck type equations
on the Wiener space}. Infin. Dimens. Anal. Quantum Probab. Relat.
Top. 13 (2010), no. 2, 273--304.

\bibitem{Malliavin97} P. Malliavin, {\it Stochastic analysis}, Grund.
Math. Wissen., vol. 313, Springer, 1997.

\bibitem{Ni04a} Lei Ni, {\it The entropy formula for linear heat
equation}. J. Geom. Anal. 14 (2004), no. 1, 87--100.

\bibitem{Ni04b} Lei Ni, {\it Addenda to ``The entropy formula for linear heat
equation''}. J. Geom. Anal. 14 (2004), no. 2, 369--374.

\bibitem{Nualart} D. Nualart, {\it The Malliavin calculus and related
topics}. Second edition. Probability and its Applications (New
York). Springer-Verlag, Berlin, 2006.

\bibitem{Shigekawa} I. Shigekawa, {\it The de Rham-Hodge-Kodaira decomposition
on abstract Wiener space}. J. Math. Kyoto Univ. 26 (1986), 191--202.

\bibitem{Sturm1} K.Th. Sturm, {\it On the geometry of measures spaces. I}. Acta Math. 196 (2006),
65--131.

\bibitem{Sturm2} K.Th. Sturm, {\it On the geometry of measures spaces. II}. Acta Math. 196 (2006),
133--177.

\bibitem{Sugita} H. Sugita, {\it Positive generalized Wiener functionals and
potential theory over abstract Wiener space}. Osaka J. Math. 25
(1988), 665--696.

\end{thebibliography}
\end{document}